\documentstyle [12pt,twoside, epsf]{article}

\def\picill#1by#2(#3)
{\vbox to #2
{\hrule width #1 height 0pt depth 0pt
\vfill\epsffile{#3}}}

\begin{document}

\date{}

\title{\Large\bf Review of ``Knots" by Alexei Sossinsky, Harvard University Press, 2002, ISBN 0-674-00944-4}

\author{Louis H. Kauffman \\ Department of Mathematics, Statistics and Computer
Science \\ University of Illinois at Chicago \\ 851 South Morgan Street\\
Chicago, IL, 60607-7045}

\maketitle

\thispagestyle{empty}

\section{Introduction}
This is a brilliant and sharply written little book about knots and theories of knots.
Listen to the author's preface.
\bigbreak

"Butterfly knot, clove hitch knot, Gordian knot, hangman's knot, vipers' tangle
- knots are familiar objects, symbols of complexity, occasionally metaphors for 
evil. For reasons I do not entirely understand, they were long ignored by
mathematicians. A tentative effort by Alexandre-Th[e]ophile Vandermonde at the end
of the eighteenth century was short-lived, and a preliminary study by the young
Karl Friedrich Gauss was no more successful. Only in the twentieth century did
mathematicians apply themselves seriously to the study of knots. But until the 
mid 1980s, knot theory was regarded as just one of the branches of topology: 
important, of course, but not very interesting to anyone outside a small
circle of specialists (particularly Germans and Americans).
\bigbreak

Today, all that has changed. Knots - or more accurately, mathematical theories of
knots - concern biologists, chemists, and physicists. Knots are trendy.
The French ``nouveaux philosophes" (not so new anymore) and postmodernists
even talk about knots on television, with their typical nerve and incompetence.
The expressions ``quantum group" and ``knot polynomial" are used indiscriminately
by people with little scientific expertise. Why the interest? Is it a passing 
fancy or the provocative beginning of a theory as important as relativity or
quantum physics?"
\bigbreak

Sossinsky continues for 119 pages in that quick, vivid, irreverent way, 
telling the story of knots and knot theory from practical knot tying to 
speculations about the relationship of knots and physics. I will first give a personal sketch of 
some history and ideas  of knot theory. Then we will return to a discussion of Sossinsky's remarkable 
book. 
\bigbreak

\section{ A Sketch of Knot Theory}
Knot theory had its most recent 
beginnings in the nineteenth century due to the curiosity of Karl Friedrich Gauss,
James Clerk Maxwell and Peter Guthrie
Tait, and the energy of Lord Kelvin (Sir William Thomson). The latter had the popular physical
theory of the day with the hypothesis that atoms were three dimensional knotted vortices in the 
all pervading ether of space. Thomson enlisted the aid of mathematicians Tait, Little and Kirkman
to produce the first tables of knots, with the hopes that these tables would shed light on the
structure of the chemical elements. Eventually, this theory foundered, first on the vast prolixity 
of knotted forms and later in the demise of the etheric point of view about the nature of space.
But at the same time, mathematical concepts of manifolds and topology were coming forth in the 
hands of Gauss, Riemann and later Poincare. With these tools it became possible to analyze 
topological phenomena such as knotting and the properties of three dimensional manifolds. It was 
Poincare's fundamental group (of a topological space) that became the first significant tool in 
knot theory. From the properties of the fundamental group of the complement of the knot, Max Dehn
(in the early 1900's) was able to prove the knottedness of the trefoil knot and its ineqivalence to 
its mirror image. In this way the deep question of detecting the chirality of knots was born.
In the 1920's James W. Alexander of Princeton University discovered a polynomial invariant of knots 
and links \cite{Alexander} that enabled many extensive computations. Alexander's polynomial could not distinguish 
knots from their mirror images, but it was remarkably effective in other ways, and it was based upon 
ideas from the fundamental group and from the structure of covering spaces of the knot complement.
In fact, Alexander based the theory of his polynomial invariant of knots on the newly discovered
Reidemeister moves, expressing the topological equivalence relation for isotopy of knots and links 
in terms of a language of graphical diagrams. (Reidemeister's exposition of the moves can be found in his
book \cite{Reidemeister}.) Alexander's version of his polynomial was expressed by the  determinant of a matrix that one
can read directly from the knot diagram. Invariance of the polynomial is proved by examining how this determinant behaved
under the Reidemeister moves. Reidemeister wrote the first book on knot theory and based it upon his 
moves.
\bigbreak

It would take another sixty years to realize the 
power inherent in Reidemeister's approach to knot theory. Topology evolved from the 1920's onward 
with the seminal work of Seifert \cite{Seifert} on knots and three manifolds and the rapid evolution of 
algebraic topology. In the early 1950's the precise role of the fundamental group in knot theory 
was made clear by the work of Ralph Fox \cite{Fox} and his students. Fox showed how one could extract the 
analogues of Alexander polynomials directly from the presentations of the groups by a remarkable 
algebraic technique (derived from the theory of covering spaces) called the free differential 
calculus. It was a non-commutative and discrete version of Newtonian calculus, adapted to the needs 
of algebraic topology and combinatorial group theory. 
\bigbreak

Then in the late 1960's John Horton Conway
published a startling paper \cite{Conway} in which he showed how to compute Alexander polynomials without any 
matrices, free calculi or determinants. His method relied on a recursive formula that expressed the
Conway version of the Alexander polynomial in terms of simpler knots and links. The Conway skein theory was met by
puzzlement on the part of topologists. It took about ten years for knot theorists to start thinking about the Conway
approach, and the first thoughts became proofs of  various sorts that the Conway method was valid and that it resulted
in a normalized version of  the Alexander polynomial. The author of this review was one of those captured souls, who
puzzled  about the Conway magic. He found two approaches to it. The first \cite{Kauffman1} went back to techniques of
Seifert. The second \cite{Kauffman2} went back to Alexander's original paper. The second approach rewrote and normalized
Alexander's determinant, converting it to a sum over combinatorial states of the link diagram. This sum over states made
it easy to prove that the resulting polynomials satisfy Conway's identities and that they are invariant under the
Reidemeister moves. The state sum produced a  fully combinatorial (graph theoretic) way to understand Alexander's
original determinant. The state summation is analogous to certain sums in graph theory and to partition functions in
statistical mechanics.  The full significance of this analogy was not apparent in 1980/81 when these relations were
discovered.
\bigbreak

At this point it is worth making a digression about the Reidemeister moves. In the 1920's Kurt Reidemeister
proved an elementary and important theorem that translated the problem of determining the topological type of a knot
or link to a problem in combinatorics. Reidemeister observed that any knot or link could be represented by a 
{\em diagram} where a diagram is a graph in the plane with four edges locally incident to each node, and with extra
structure at each node that indicates an over-crossing of one local arc (consisting in two local edges in the graph)
with another. See Figure 1. The diagram of a classical knot or link has the appearance of a sketch of the knot, but it
is a rigorous and exact notation that represents the topological type of the knot. Reidemeister showed 
that two diagrams represent the same topological type (of knottedness or linkedness) if and only if one diagram can be
obtained from another by planar homeomorphisms coupled with a finite sequence of the {\em Reidemeister moves} \
illustrated in Figure 2. Each of the Reidemeister moves is a local change in the diagram that is applied as shown in 
this Figure.

$$ \picill3inby2in(Knot)  $$

\begin{center}
{\bf Figure 1 - A Knot Diagram}
\end{center}

$$ \picill3inby3in(Moves)  $$

\begin{center}
{\bf Figure 2 - Reidemeister Moves}
\end{center}

The first move is special for a number of reasons. One can permute the performance of the first move with the other 
moves. So one can save up the doings of the first Reidemeister move until the very end of a process, if one so 
desires. And from a physical point of view, there are good reasons to not use the first Reidemeister move. The move
is designed to remove a curl in the line, but a curl in a rope does not go away, it just gets hidden as a twist when 
you pull on the rope. The reader should try this with a bit of string or a rubber band. Form a curl as in Figure 3 and 
then pull on the ends of the string or band. You will find that the curl is transmuted into a twist, and if you relax the
string or band, the curl can reappear. For this reason, it is useful to consider just the equivalence relation
generated by the second and third Reidemeister moves. This relation is called {\em regular isotopy}.
Regular isotopy was first defined by Bruce Trace in \cite{Trace} and has turned out to be a useful companion to the 
full equivalence relation defined by all three Reidemeister moves. One way of thinking about regular isotopy is that 
one is talking about a {\em framed knot or link} where, by framing one means that there is an embedding of a band
(i.e. the cross product of a circle with the unit interval) for each component of the link, so that the bands do not
intersect one another. The bands can twist and this twisting models the twisting of a physical rope or a rubber band.
In order to model embedded bands with curly diagrams we actually do not just remove the first Reidemeister move, but
we replace it by the move illustrated in Figure 3. Note that there are two ways, shown in this Figure, to make a
curl out of a full twist in a band. These two curls are equated with one another and this becomes the replacement for
the first Reidemeister move. With this replacement, one refers to the equivalence classes as framed links in the 
blackboard framing (since these diagrams can be drawn on a blackboard). We shall speak about ``measuring the framing"
when speaking about the curls in a diagram taken up to regular isotopy because of this connection with framed links.
\bigbreak

$$ \picill3inby4in(Framing)  $$

\begin{center}
{\bf Figure 3 - Framing Equivalence}
\end{center}

In 1984 Vaughan Jones dropped a bombshell \cite{Jones} from which knot theory and indeed modern mathematics has
yet to recover. By following an analogy between the structure of Artin Braid Group and certain
identities in a class of von Neumann algebras, Jones discovered new representations of the braid group
and used these representations to produce an entirely new Laurent polynomial invariant of knots
and links. On top of this, Jones showed that his new (one variable) invariant satisfied a skein
relation that was almost the same as the relation for the Conway polynomial. Only the coefficients
were changed. This was shocking. On top of that, the Jones polynomial could distinguish many knots 
from their mirror images, leaving the Alexander polynomial in the dust.
\bigbreak

Jones' invariant was quickly generalized to a two-variable polynomial invariant of knots and links that goes by the 
acronym ``Homflypt" polynomial after the people who discovered it: Hoste, Ocneanu, Millett, Freyd, Lickorish, Yetter,
Przytycki and Traczyk. There were collaborations with Millett and Lickorish working together, Freyd and Yetter working
together and Przytycki and Traczyk working together. This generalization was proved to have its properties in a number
of different ways: direct induction on the properties of knot and link diagrams, algebras related to skeins of knots
and links, representations of Hecke algebras (generalizing the von Neumann algebras used by Jones), some category theory
in the bargain. A few months went by, and then Brandt, Millett and Ho discovered a new one-variable invariant of
unoriented knots and links with a different skein relation. This invariant did not detect the difference between knots
and their mirror images, but there was a new idea in it, namely that one could smooth an unoriented crossing in two
different ways (as shown in Figure 4). When the author of this review received a copy of the announcement of this
invariant, he was inspired and astonished to realize that their invariant could be 
generalized to a two variable polynomial invariant $L_{K}(z,a)$ of knots and links that {\em did detect} the difference
between  many links and their mirror images.  The key idea for this generalization is to regard the original invariant
as an  invariant of {\em framed links}, adding an extra variable that measures the {\em framing}. The reviewer had
earlier found an analogous way to understand the Homflypt polynomial. This was not published until later in
\cite{Kauffman3, Kauffman4}.

$$ \picill4inby2.5in(Smoothings) $$
\begin{center} {\bf Figure 4 - Smoothings} 
\end{center}

The invariant $L_{K}$ satisfies the following formulas 

$$L_{\mbox{\large $\chi$}} + L_{\overline{\mbox{\large $\chi$}}} = z( L_{\mbox{\large $\asymp$}} +  L_{)(})$$
$$L_{\mbox{\large $\gamma$}} = a L_{\smile} \hspace{.5in}$$ and
$$L_{\overline{\mbox{\large $\gamma$}}} = a^{-1} L_{\smile} \hspace{.5in}$$
where the small diagrams represent parts of larger diagrams that are identical except  at
the site indicated in the diagram. We take the convention that the letter chi, \mbox{\large $\chi$},
denotes a crossing where {\em the curved line is crossing over the straight
segment}. The barred letter denotes the switch of this crossing, where {\em the curved
line is undercrossing the straight segment}. 
In this formula we have used the notations $\mbox{\large $\asymp$}$ and $)($ to indicate the 
two new diagrams created by the two smoothings of a single crossing in the diagram $K$. That is, the 
four diagrams differ at the site
of one crossing in the diagram $K$.  
Here \mbox{\large $\gamma$}  denotes a curl of positive type as indicated in Figure 5, 
and  $\overline{\mbox{\large $\gamma$}}$ indicates a curl of negative type, as also seen in this
figure. The type of a curl is the sign of the crossing when we orient it locally. Our convention of
signs is also given in Figure 5. Note that the type of a curl  does not depend on the orientation
we choose.  The small arcs on the right hand side of these formulas indicate
the removal of the curl from the corresponding diagram.

$$ \picill4inby2.5in(Curls) $$
\begin{center} {\bf Figure 5 - Crossing Signs and Curls} 
\end{center}
\vspace{3mm}

\noindent The polynomial $L_{K}(z,a)$ is invariant under regular isotopy and can be  normalized to an invariant of
ambient isotopy by the definition  $$F_{K}(z,a) = a^{-w(K)}L_{K}(z,a),$$ where we chose an orientation for $K$, and
where $w(K)$ is  the sum of the crossing signs  of the oriented link $K$. $w(K)$ is called the {\em writhe} of $K$. 
The convention for crossing signs is shown in  Figure 5.
\bigbreak

This use of regular isotopy is a key ingredient in defining the polynomial $L_{K},$ and it turns out to be important
in defining many other knot polynomials as well. In the case of the reviewer's research, discovering $L_{K}$ 
opened a doorway to finding a remarkably simple model of the original Jones polynomial the {\em bracket state sum
model} \cite{Kauffman5, Kauffman6}. The bracket state sum is also an invariant of regular isotopy, and can be normalized
just like
$L_{K}$ by using the writhe. The idea behind the bracket state sum is the notion that one might form an invariant of
knots and links  by summing over ``states" of the link diagram in analogy with summations over states of physical
systems (called  partition functions) that occur in statistical mechanics. The reviewer had earlier discovered a model
for the Alexander- Conway polynomial in this form \cite{Kauffman2} and was convinced that such models should exist for
the new knot polynomials. The first case of such a model occurs with the bracket state summation. 
\bigbreak

The {\em bracket polynomial} , $<K> \, = \, <K>(A)$,  assigns to each unoriented link diagram $K$ a 
Laurent polynomial in the variable $A$, such that
   
\begin{enumerate}
\item If $K$ and $K'$ are regularly isotopic diagrams, then  $<K> \, = \, <K'>$.
  
\item If  $K \amalg O$  denotes the disjoint union of $K$ with an extra unknotted and unlinked 
component $O,$ then 

$$< K \amalg O> \, = \delta<K>,$$ 
where  $$\delta = -A^{2} - A^{-2}.$$
  
\item $<K>$ satisfies the following formulas 

$$<\mbox{\large $\chi$}> \, = A <\mbox{\large $\asymp$}> + A^{-1} <)(>$$
$$<\overline{\mbox{\large $\chi$}}> \, = A^{-1} <\mbox{\large $\asymp$}> + A <)(>,$$
\end{enumerate}

\noindent where the small diagrams represent parts of larger diagrams that are identical except  at
the site indicated in the bracket. We take the same conventions as described for the $L-$ polynomial.
Note, in fact that it follows from the bracket formulas that 
$$ <\mbox{\large $\chi$}> \, + <\overline{\mbox{\large $\chi$}}> \, = (A+A^{-1})(<\mbox{\large $\asymp$}> + <)(>).$$
Making the bracket polynomial a special case of the $L-$ polynomial.
\bigbreak

\noindent It is easy to see that Properties $2$ and $3$ define the calculation of the bracket on
arbitrary link diagrams. The choices of coefficients ($A$ and $A^{-1}$) and the value of $\delta$
make the bracket invariant under the Reidemeister moves II and III (see \cite{Kauffman5}). Thus
Property $1$ is a consequence of the other two properties. 

\bigbreak

In order to obtain a closed formula for the bracket, we now describe it as a state summation.
Let $K$ be any unoriented link diagram. Define a {\em state}, $S$, of $K$  to be a choice of
smoothing for each  crossing of $K.$ There are two choices for smoothing a given  crossing, and
thus there are $2^{N}$ states of a diagram with $N$ crossings.
 In  a state we label each smoothing with $A$ or $A^{-1}$ according to the left-right convention 
discussed in Property $3$ (see Figure 4). The label is called a {\em vertex weight} of the state.
There are two evaluations related to a state. The first one is the product of the vertex weights,
denoted  

$$<K|S>.$$
The second evaluation is the number of loops in the state $S$, denoted  $$||S||.$$
  
\noindent Define the {\em state summation}, $<K>$, by the formula 

$$<K> \, = \sum_{S} <K|S>\delta^{||S||-1}.$$
It follows from this definition that $<K>$ satisfies the equations
  
$$<\mbox{\large $\chi$}> \, = A <\mbox{\large $\asymp$}> + A^{-1} <)(>,$$
$$<K \amalg  O> \, = \delta<K>,$$
$$<O> \, =1.$$
  
\noindent The first equation expresses the fact that the entire set of states of a given diagram is
the union, with respect to a given crossing, of those states with an $A$-type smoothing and those
 with an $A^{-1}$-type smoothing at that crossing. The second and the third equation
are clear from the formula defining the state summation. Hence this state summation produces the
bracket polynomial as we have described it at the beginning of the  section. 

\bigbreak

In computing the bracket, one finds the following behaviour under Reidemeister move I: 
  $$<\mbox{\large $\gamma$}> = -A^{3}<\smile> \hspace {.5in}$$ and 
  $$<\overline{\mbox{\large $\gamma$}}> = -A^{-3}<\smile> \hspace {.5in}$$

\noindent where \mbox{\large $\gamma$}  denotes a curl of positive type as indicated in Figure 5, 
and  $\overline{\mbox{\large $\gamma$}}$ indicates a curl of negative type, as also seen in this
figure. The type of a curl is the sign of the crossing when we orient it locally. Our convention of
signs is also given in Figure 5. Note that the type of a curl  does not depend on the orientation
we choose.  The small arcs on the right hand side of these formulas indicate
the removal of the curl from the corresponding diagram.  

\bigbreak
  
\noindent The bracket is invariant under regular isotopy and can be  normalized to an invariant of
ambient isotopy by the definition  
$$f_{K}(A) = (-A^{3})^{-w(K)}<K>(A),$$ where we chose an orientation for $K$, and where $w(K)$ is 
the {\em writhe} of $K$.  By a change of variables one obtains the original
Jones polynomial, $V_{K}(t),$  for oriented knots and links from the normalized bracket:

$$V_{K}(t) = f_{K}(t^{-\frac{1}{4}}).$$

\noindent The bracket model for the Jones polynomial is quite useful both theoretically and in terms
 of practical computations. One of the neatest applications is to simply compute $f_{K}(A)$ for the
trefoil knot $T$ and determine that  $f_{T}(A)$ is not equal to $f_{T}(A^{-1}) = f_{-T}(A).$  This
shows that the trefoil is not ambient isotopic to its mirror image, a fact that is quite tricky to
prove by classical methods. To this day, it is still an open problem whether there are any non-trivial classical knots
with unit Jones polynomial.
\bigbreak

After the bracket polynomial model demonstrated the idea of a state summation directly related to the knot diagram
there came lots of other examples of state summations, using algebraic machinery from statistical mechanics and 
Hopf algebras (See \cite{Kauffman7} for an exposition of some of this development). The subject grew rapidly and then
underwent another change in the late 1980's when Edward Witten showed how to think about such invariants in terms of
quantum field theory
\cite{Witten}. Witten's methods indicated that one should be able to construct invariants of three-manifolds, and in
fact David Yetter \cite{Yetter1, Yetter2} had shown a similar pattern earlier using categories and Reshetikhin and Turaev
\cite{RT} showed explicitly how to accomplish this goal using the algebra of quantum groups (Hopf algebras).  Witten's
approach brought gauge field theory into the subject of knot invariants. In Witten's approach there is  a Lie algebra
valued connection (a differential one-form) defined on three dimensional space. One can measure the  trace of the
holonomy of this connection around a knotted loop $K$. This measurement is called the {\em Wilson loop}
$W_{K}(A).$ The Wilson loop is not itself a topological invariant of $K$, but a suitable average of $W_{K}(A)$ over
many connections $A$ should be a (framed) invariant. Witten suggested the specific form of this averaging process, and
indeed that average works just as advertised at a formal level. The formalism leads to the idea of a topological
quantum field theory  due jointly to Witten and Atiyah \cite{Atiyah}, a concept that has been highly influential since
that time. Witten's techniques
are highly suggestive of fully rigorous combinatorial approaches. His work
continues to act as a catalyst for new approaches to the structure of the invariants, and as a connection of the subject
to mathematical physics. One of the most exciting aspects of the Wilson loop formulation is that it leads to connections
between knot theory and theories of quantum gravity \cite{Smolin1, Smolin2} and string theory \cite{Ooguri, LM}.
\bigbreak

Other approaches to link invariants grew out of Witten's work, most notably the theory of Vassiliev 
invariants \cite{Vassiliev}. Vassiliev invariants were first seen as coming from considerations of the topology of the
space of all  knots and singular knots, but were then quickly connected to combinatorial topology \cite{Birman} on the
one hand, and to  Witten's work \cite{Bar-Natan1, Bar-Natan2} on the other. In  \cite{Bar-Natan2} there is an 
excellent account of the approach to Vassiliev invariants via the Kontesevich integrals and in \cite{Kauffman7,
Kauffman8} the reader will find a heuristic account telling how the Kontsevich integrals arise in the Feynman diagram
expansion of Witten's functional integral. 
The result of this evolution  has been a very clear view of just how it is that Lie algebraic structures are related to
invariants of knots and links. The reader should recall that a Lie algebra is a linear algebra with a non-associative
multiplication (here denoted $ab$) such that
$ab=-ba$ for all
$a$ and
$b$ in the algebra, and such that $a(bc) = (ab)c + b(ac)$ (the Jacobi identity) for all elements $a,b,c$ in the algebra.
A diagrammatic version of the Jacobi identity can be seen just beneath the surface of the Reidemeister moves, when one
takes the Vassiliev point of view, and it is this transition through diagrammatic (or categorical) algebra that makes a
deep connection between knot theory, algebra and mathematical physics. 
\bigbreak

\bigbreak Three recent developments are worth mentioning. First is the discovery of a generalization (categorification)
of the bracket  polynomial by Misha Khovanov \cite{Khovanov, Bar-Natan3} that writes the original Jones polynomial as an
Euler characteristic of a complex whose graded cohomology leads to new invariants and to invariants of surfaces in four-
dimensional space. One extraordinary outgrowth of the Khovanov work is the deep work of Ozwath and Szabo \cite{OZSA}
finding a categorification of  the Alexander polynomial that is related to the state sum in \cite{Kauffman2}. The
second is the discovery by Morwen Thistlethwaite
\cite{Morwen} of examples of non-trivial links that have the same Jones polynomial as trivial links, and the extension
of this result by Eliahou, Kauffman and Thistlethwaite \cite{EKT} to infinite families of links with this property. The
third is the discovery of virtual knot theory \cite{Kauffman9}, a generalization of classical knot theory to knots in
abstract surfaces where there exist inifinitely many non-trivial virtual knots with unit Jones polynomial. We do not yet
know if any of these virtual counterexamples will yield classical knots with unit Jones polynomial. The fundamental
problem remains open.

\section{Sossinsky's ``Knots"}
The book ``Knots"  \cite{Sossinsky} by Alexei Sossinsky is for general readers. It is a companion
piece for the book ``Knots, Links, Braids and $3$-Manifolds - An Introduction to
the New Invariants in Low-Dimensional Topology" \cite{PS} by Prasolov and Sossinsky, published in translation by
the American Mathematical Society. The latter book is not under review here, but is recommended to all readers who 
enjoyed reading the book that is under review.
Sossinsky's ``Knots" begins with a very nicely illustrated tour of practical knots and decorative knotwork. It is
recalled how Lord Kelvin (Sir William Thomson) in the 1860's had a theory of vortex atoms in which the atomic 
constituents of matter were to be modeled on three dimensional vortices, knots in the ether. This theory 
inspired Kelvin to enlist the aid of mathematicians Tait, Kirkman and Little to construct the first tables of 
knots. The vortex theory eventually came to a stop with the rejection of the ether in favor of relativity, but the
idea of knots related to physics lives on to this day. Sossinsky's book describes this beginning for knot theory
and then continues with a discussion of topological equivalence, wild knots, braids and the braid group and a proof of
the basic theorem of Alexander that any knot can be represented as a closed braid. Sossinsky then goes on to describe
a modern algorithm due to Pierre Vogel for transforming a knot into a braid, and a new algorithm of Dehornoy for finding 
a minimal representation of a braid. Chapter 3 discusses planar diagrams of knots and the Reidemeister moves. Chapter 4
discusses the arithmetic of knots under connected sum, and gives a sketch of one proof of the theorem that you cannot
cancel knots. That is, if $K$ and $K'$ have an unknotted connected sum, then $K$ and $K'$ are individually unknotted.
This is a fact that people often find surprising at first and the proof given by Sossinsky uses a technique known to 
topologists as ``the method of infinite repetition". This method of proof is probably more astonishing than the Theorem!
The chapter ends with a discussion of the factorization of knots into prime knots. Here is the beginning of a fruitful
analogy of knots and numbers that is really just beginning to be explored. 
\bigbreak

Chapter 5 is packed with a number of things.
There is a discussion of Conway's approach and reformulation of the Alexander polynomial. Conway's approach uses 
oriented diagrams and involves two operations that can be performed at a diagrammatic crossing. One can switch the two
strands, or one can smooth the crossing by reattaching the strands, preserving orientation and removing the 
crossover in the process. These two operations are reminiscent of the operations of topoisomerase enzymes and of 
combinatorics of DNA recombination (as was observed some years after Conway devised his methods). Sossinsky takes the
opportunity to make a quick digression into the subject of DNA topology. He then discusses how one calculates with 
Conway's polynomial, and how natural it is to generalize (with 20/20 hindsight!) this construction of Conway to the 
2-variable Homflypt polynomial. Chapter 6 returns the discussion to the original Jones polynomial and constructs the 
bracket polynomial state summation in some detail, including its relationship with the Potts model in statistical 
mechanics. The reviewer is pleased with this exposition of the bracket. The reader will find some speculation by 
Sossinsky on how this construction was discovered. In this review the author of the review has given his version of this
story in section 2. Chapter 7 goes on to discuss the Vassiliev invariants. The discussion of this advanced topic is
accomplished neatly, with an emphasis on the ideas and the diagrammatic context. In fact, Sossinsky discusses the
so-called four term relation for the Vassiliev invariants, and in this way comes right up to the fundamental
relationship between knot invariants and Lie algebras. This is an example of the daring nature of this popular
exposition. Chapter 8 treats relationships of knot theory with physics, including a bit about statistical mechanics,
quantum field theory and quantum gravity.
\bigbreak

This book is a tour de force. It is a great read and, I believe that it is a book that can actually be understood by 
a very wide readership. When large scale topics are discussed, Sossinsky deftly, sometimes ironically, places key ideas
before the reader. A myriad of details will have to be dealt with by one who wants to master these topics but here the
ideas are laid bare. I would recommend this book as a first book on knot theory to anyone who asks.
\bigbreak

There are a few slips and misprints. I will tell them here to the best of my knowledge. This is a book that deserves 
many editions and it is to be assumed that the author will correct these few stumbling blocks in the very next edition 
of the book. Page 43, Figure 3.4 shows a diagram that actually can be reduced to the unknot by Reidemeister moves
without making the number of crossings increase (contrary to the claim in the book). Diagrams with that property are 
not hard to construct, but this is not one of them. On page 70 in Figure 5.7, $A$ is a figure eight knot, contrary to
the book's assertion, and $B$ is a trefoil knot. On page 71 it is claimed that the Conway polynomial of the  trefoil knot
and the figure eight knot are the same. This is not the case. In the book's notation we have
$1+x^2$ for the Conway polynomial of the trefoil knot and $1-x^2$ for the Conway polynomial of the figure eight
knot. It is important for the next edition to correct this error, since a beginner in calculating these polynomials
could become mighty confused by such a claim. On page 73 it is said that Louis Kauffman is at the University of Chicago.
In fact, he is at the University of Illinois at Chicago.
\bigbreak

I cannot resist ending this review with one more quote from the book. Sossinsky expresses throughout  his 
amazing enthusiasm for mathematics. In this quote he is speaking of the last touches, completing
the construction and proof of invariance of the bracket polynomial. He says 
\bigbreak

`` God knows I do not like exclamation points. I generally prefer Anglo-Saxon understatement to the exalted
declarations of the Slavic soul. Yet I had to restrain myself from putting two exclamation points instead of just
one at the end of the previous section. Why? Lovers of mathematics will understand. For everyone else: the emotion
a mathematician experiences when he encounters (or discovers) something similar is close to what the art lover
feels when he first looks at Michelangelo's {\em Creation} in the Sistine Chapel. Or better yet (in the case of 
a discovery), the euphoria that the conductor must experience when all the musicians and the choir, in the same
breath that he instills and controls, repeat the ``Ode to Joy" at the end of the fourth movement of Beethoven's 
{\em Ninth}." 
\bigbreak

\noindent I could not agree more. Read the book!
 \bigbreak

\end{document}